\documentclass[12pt,a4paper]{amsart}
\usepackage{amsfonts,amsmath,amssymb}

\def\GL{\operatorname{\mathrm {GL}}\nolimits}
\def\SO{\operatorname{\mathrm {SO}}\nolimits}

\def\SL{\operatorname{\mathrm {SL}}\nolimits}
\def\sl{\operatorname{\mathfrak {sl}}\nolimits}
\def\Lie{\operatorname{\mathrm {Lie}}\nolimits}
\def\Aut{\operatorname{\mathrm {Aut}}\nolimits}
\def\Rad{\operatorname{\mathrm {Rad}}\nolimits}
\def\diag{\operatorname{\mathrm {diag}}\nolimits}
\def\sign{\operatorname{\mathrm {sign}}\nolimits}
\def\sgn{\operatorname{\mathrm {sign}}\nolimits}
\def\vol{\operatorname{\mathrm {vol}}\nolimits}
\def\Hei{\operatorname{\mathrm {Heis}}\nolimits}
\def\Heis{\operatorname{\mathrm {Heis}}\nolimits}
\def\Cent{\operatorname{\mathrm {Cent}}\nolimits}
\def\ind{\operatorname{\mathrm {Ind}}\nolimits}
\def\tr{\operatorname{\mathrm {tr}}\nolimits}
\def\ad{\operatorname{\mathrm {ad}}\nolimits}

\def\bbR{\mathbb R}
\def\bbC{\mathbb C}
\def\bbS{\mathbb S}
\def\bbT{\mathbb T}

\newtheorem{thm}{Theorem}[section]
\newtheorem{prop}[thm]{Proposition}
\newtheorem{lem}[thm]{Lemma}

\def\beginpf{\textsc{Proof.}}
\def\endpf{\hfill$\Box$}
\title[Poisson Summation and Endoscopy for $\SL(3,\mathbb R)$]{Poisson Summation and Endoscopy\\ for $\SL(3,\mathbb R)$}
\author[Do Ngoc Diep, Do Thi Phuong Quynh]{Do Ngoc Diep${}^{1}$, Do Thi Phuong Quynh${}^{2}$}
\address{${}^{1}$ Institute of Mathematics, VAST, 18 Hoang Quoc Viet road, Cau Giay district, 10307 Hanoi, Vietnam \newline
{\tt Email: dndiep@math.ac.vn}}
\address{${}^{2}$ {\sc Medicine University Colllege, Thai Nguyen Universtiy, Thai Nguyen City, Vietnam}\newline
{\tt Email: phuongquynhtn@gmail.com}}
\begin{document}
\date{\bf Version of \today}
\maketitle
\begin{abstract} The group is interesting as the first example of split rank 2 semisimple group, all the irreducible unitary representations of which are known. 
We make a precise realization of the discrete series representations (in Section 2) by using the Orbit Method and Geometric Quantization, a computation of their traces (Section 3) and an exact formula for the noncommutative Poisson summation and endoscopy of  for this group (in Section 4).

\textit{Key terms}: trace formula; orbital integral; transfer; endoscopy

\textbf{2010 Mathematics Subject Classification}: \textit{Primary}: 22E45; \textit{Secondary}: 11F70, 11F72, 
\end{abstract}
\tableofcontents
\section{Introduction}
The questions of finding all the irreducible unitary representations of a reductive Lie group and decomposition of a particular represention into a discrete or continuous direct sum or integral of irreducible ones  are the basic questions of harmonic analysis on reductive Lie groups. In particlular, the discrete part of the regular representation of reductive Lie groups is the discrete sum of discrete series representations. 

Often these problems are reduced to the trace formula, because, as known, the unitary representations are uniquely defined by it generalized character and infinitesimal character. Following Harish-Chandra, the generalized character is defined by its restriction to the maximal compact subgroup, as the initial eigenvalue problem for the generalized Laplacian (the Casimir operator) with the infinitesimal character as the eigenvalues infinitesimal action of Casimir operators.

There is a very highly developed theory of Arthur-Selberg trace formula. The theory is complicated and one reduces it to the same problem for smaller   endoscopic subgroups. It is called the transfer and plays a very important role in the theory.
By definition, an endoscopic subgroup is the connected component of the centralizer of regular semisimple elements, associated to representations, namely by the orbit method.

For the discrete series representations of $\SL(3,\mathbb R)$ the following endoscopic groups should be considered:
\begin{itemize}
\item \textit{the elliptic case}: diagonal subgroup of regular elements
$$H= \left\{\diag(a_1,a_2,a_3); \quad a_1a_2a_3 = 1, a_i\ne a_j,\forall i \ne j\right\}\cong \bbT^2;$$
\item \textit{the parabolic case}: blog-diagonal subgroup of regular elements 
$H = \left\{\begin{pmatrix} \cos\theta & \sin\theta & 0\\ -\sin\theta & \cos\theta & 0 \\ 0 & 0 & 1\end{pmatrix}\right\}\cong \SO(2)$;
\item the trival case: the group $H= \SL(3,\bbR)$ its-self.
\end{itemize}
We show in this paper that the method that J.-P. Labesse \cite{labesse} used for $\SL(2,\bbR)$ is applied also for $\SL(3,\bbR)$.
Therefore one deduces the transfer formula for the discrete series representations and limits of  $\SL(3,\bbR)$ to the corresponding endoscopic group. 

For the group $\SL(3,\bbR)$ we make a precise realization of the discrete series representations (in Section 2) by using the Orbit Method and Geometric Quantization to the solvable radical, a computation in the context of $\SL(3,\bbR)$ of their traces (Section 3) and an exact formula for the noncommutative Poisson summation and endoscopy of  for this group (in Section 4).

\section{Irreducible Unitary Representations of $\SL(3,\mathbb R)$}

\subsection{The structure of $\SL(3,\mathbb R)$} The following notions and results are folklore and we recall them to fix an appropriate system of notations.

Let us remind that the group  $\SL(3,\mathbb R)$ is 
$$\SL(3,\mathbb R= \left\{ X \in \GL(3,\mathbb R) | \det X = 1  \right\}.$$ 

Denote by $\mathfrak{sl}(3,\mathbb R)$ its Lie algebra $\Lie\SL(3,\mathbb R)$, $\theta$ the Cartan involution of the group $G= \SL(3,\mathbb R)$ which is $\theta(X) = {}^tX^{-1}$ . The corresponding Cartan involution of   for its Lie algebra $\mathfrak{sl}(3,\mathbb R)$ is denoted by the same symbol $\theta\in \Aut\mathfrak{sl}(3,\mathbb R)$,  $$\theta(X) = -{}^tX, X\in \sl(3,\mathbb R)$$ The maximal compact subgroup $K$ of $G$ is the orthgonal group
$$K = \SO(3)$$ 
 is the subgroup of $G$, the Lie algebra $\mathfrak k$ of which is consisting of all the matrices on which the Cartan involution has eigenvalue $+1$,
$$\mathfrak k = \left\{\left. X\right| \theta(X) = -{}^t X = X    \right\}.$$ 
The Borel subgroup of $\SL(3,\mathbb R)$ is the minimal parabolic subgroup $$P_0 = B=\left\{\left. p= \begin{pmatrix} m & *\\ 0 & \det m^{-1} \end{pmatrix}\right| m \in U(2) \right\},$$  
the Lie algebra of which is consisting of all the matrices with eigenvalue $-1$,
$$\mathfrak b = \left\{\left.  X\in \mathfrak g \right| \theta(X) = -X   \right\}$$. Up to conjugacy, there are two possible Borel subgroups: the split Borel subgroup 
$$B_s=\left\{\left.\begin{pmatrix} t_1 & * & *\\  0 & t_2 & *\\ 0 & 0 & t_3   \end{pmatrix}\right| t_i \in \mathbb R_{+}, t_1t_2t_3 = 1 \right\}$$ with the maximal abelian subgroup $A= \diag(t_1,t_2,t_3)$ and unipotent radical
$$U=\left\{\begin{pmatrix} 1 & * & *\\  0 & 1 & *\\ 0 & 0 & 1   \end{pmatrix}\right\},$$
the compact subgroup of of $B$ is $$K_s = B_s \cap K =\left\{\begin{pmatrix} \pm 1 & 0 & 0\\ 0 & \pm 1 & 0\\ 0 & 0 & 1 \end{pmatrix}  \right\}\cong \mathbb Z_2 =\mathbb Z/2\mathbb Z$$
 $B_s=\mathbb Z_2AU$,
and the non-split Borel subgroup 
$$B_n=\left\{\left.\begin{pmatrix} t_1\cos \theta & t_1\sin\theta & *\\ -t_1\sin\theta & t_1\cos\theta & *\\ 0 & 0 & t_2   \end{pmatrix}\right| t_i \in \mathbb R_{+}, t_1^2t_2 = 1 \right\}$$
 has a maximal split abelian subgroup $A= \diag(t_1,t_1,t_2)$ and the unipotent radical
$$U=\left\{\begin{pmatrix} 1 & 0 & *\\  0 & 1 & *\\ 0 & 0 & 1   \end{pmatrix} \right\},$$ the maximal compact subgroup of $B_n$ is
$$K_n=K \cap B_n = \left\{\left.\begin{pmatrix} \pm\cos \theta & \pm\sin\theta & 0\\ \mp\sin\theta & \pm\cos\theta & 0\\ 0 & 0 & 1\end{pmatrix}\right| \theta\in [0,2\pi) \right\}$$
 and $B_n = K_nAU$, 

\textit{In the split case},
the group $\SL(3,\mathbb R)$ admits the well-known Cartan decomposition in form 
$G= B_sK$. The Borel subgroup  $B_s$ is endowed with a further decomposition into a semi-direct product $B=M A U$, $M=\{\pm1\}$, of a maximal split torus $A=(\mathbf R^*_+)^2$, the Lie algebra of which is
$$\mathfrak a = \left\{\left. H = \begin{pmatrix} \lambda_1 & 0 & 0\\
0 & \lambda_2& 0\\ 0 & 0 & \lambda_3\end{pmatrix} \right| \lambda_i \in \mathbb R , \lambda_1+\lambda_2+\lambda_3 =0  \right\}$$
 and the unipotent radical $U= \Rad_uB \cong \Hei_3$, generated by matrices
$$X = \begin{pmatrix} 
0 & 1 & 0\\ 
0 & 0 & 0\\
0 & 0 & 0\end{pmatrix},
Y = \begin{pmatrix} 
0 & 0 & 0\\ 
0 & 0 & 1\\
0 & 0 & 0\end{pmatrix},
Z = \begin{pmatrix} 
0 & 0 & 1\\ 
0 & 0 & 0\\
0 & 0 & 0\end{pmatrix},$$
satisfying the Heisenberg commutation relation $[X,Y]=Z$.
$$\mathfrak b = \mathfrak u \oplus \mathfrak a \oplus \mathfrak m,$$ where 
$\mathfrak u =\Lie \Heis_3 = \langle X,Y,Z\rangle$, $\mathfrak a = \langle H_1= \diag(1, -1,0), H_2=\diag(1,0,-1) \rangle$, $\mathfrak m = 0$.

\textit{In the nonsplit case}, 
the group $\SL(3,\mathbb R)$ admits the well-known Cartan decomposition in form of 
$G= B_nK$. The Borel subgroup  $B_n$ is endowed with a further decomposition into a semi-direct product $B=MAU$ of a maximal split torus $A$, the Lie algebra of which is
$$\mathfrak a = \left\{\left. H = \begin{pmatrix} \lambda_1 & 0 & 0\\
0 & \lambda_1& 0\\ 0 & 0 & \lambda_2\end{pmatrix} \right| \lambda_i \in \mathbb R , 2\lambda_1+\lambda_2 =0  \right\}$$
 and the unipotent radical $U= \Rad_uB \cong \Heis(3,\mathbb R)$, generated by matrices
$$X = \begin{pmatrix} 
0 & 1 & 0\\ 
0 & 0 & 0\\
0 & 0 & 0\end{pmatrix},
Y = \begin{pmatrix} 
0 & 0 & 0\\ 
0 & 0 & 1\\
0 & 0 & 0\end{pmatrix},
Z = \begin{pmatrix} 
0 & 0 & 1\\ 
0 & 0 & 0\\
0 & 0 & 0\end{pmatrix},$$
satisfying the Heisenberg commutation relation $[X,Y]=Z$.
$$\mathfrak b = \mathfrak u \oplus \mathfrak a \oplus \mathfrak m,$$ where 
$\mathfrak u =\Lie \Heis(3,\mathbb R) = \langle X,Y,Z\rangle$, $\mathfrak a = \langle H=\diag(1,1,-2) \rangle$, $$\mathfrak m \cap \mathfrak b= \left\langle T=\begin{pmatrix}   
0 & 1 & 0\\ -1 & 0 & 0\\ 0 & 0 & 0
\end{pmatrix}\right\rangle.$$ 
We have the commutation relations of the harmonic oscillator $$[T,X] = 2X, [T,Y]=-2Y, [T,Z] = 0$$

The center $C(\mathfrak k)$ is of dimension 1,
$$C(\mathfrak k) = \left\{\left. \begin{pmatrix}
i\lambda & 0 & 0\\
0 & i\lambda & 0\\
0 & 0 & -2i\lambda 
\end{pmatrix}  \right| \lambda \in \mathbb R, i = \sqrt{-1}  \right\}$$
There is a compact Cartan subalgebra $\mathfrak h$ consisting of all diagonal matrices $$\mathfrak h = \left\{\left.  \diag(ih_1,ih_2,ih_3 ) \right|  h_1,h_2,h_3 \in \mathbb R, h_1 + h_2 + h_3 = 0\right\}\subset \mathfrak k.$$
The associate root system is 
$$\Delta(\mathfrak g_\mathbb C , \mathfrak k_\mathbb C) = \left\{\left. \alpha_{kl} = \alpha_k - \alpha_l \right| \alpha_k(h_l) = \delta_{kl}\right\}.$$
It means that $\alpha_{kl} = (0,\dots,0,\underbrace{1}_k,0,\dots, 0,\underbrace{-1}_l,0,\dots,0 )\in \mathfrak h^* $, $1 \leq k \ne l \leq 3$.
The subroot system of compact roots is
$\Delta_c = \{\pm \alpha_{12}\} $ The noncompact root system is
$\Delta_n = \{\pm\beta,\pm 2\beta\}$  where $\beta$ is the noncompact root such that 
$$\beta( \begin{pmatrix}\lambda & 0 & 0\\
0 & \lambda& 0\\ 0 & 0 & -2\lambda\end{pmatrix}) = \lambda$$ 

The coroot system is 
$$\Delta(\mathfrak g_\mathbb C , \mathfrak k_\mathbb C)^* = \left\{ H_{kl} = E_{kk} - E_{ll} \left| {E_{kl} = \mbox{ elementary matrix with the }\atop \mbox{only nonzero entry 1 on position } (i,j)}\right. \right\}$$

\begin{thm}
The discrete series representations of $\SL(3,\mathbb R)$ is obtained $\ind_B^G(\pi^k \otimes \chi_{\lambda}^{\pm})$ by induction from $B$ to $G$ by the tensor product of an irreducible representation of highest weight $k$ irreducible representation $\pi^k$ of the normalizer  $M_K = K \cap M$ of a semisimple element of $A$ in the maximal compact subgroup $K$,  and a character  $\chi_{\lambda_0}^\pm(uak) = a^{i\lambda} (\sgn a)^\varepsilon, \varepsilon=0,1 $ of the split component $A$ of $B$,
$$\sigma^\pm_k = \pi^k \otimes \chi_{\lambda}^{\pm}.$$
\end{thm}
\textsc{Proof.}
The proof is consisting of geometric realization of the discrete series representations of $\SL(3,\mathbb R)$. Because the Borel subgroup is a minimal parabolic subgroup, 
the coadjoint action of $K$ in $\mathfrak g^*$ keeps the set $\mathfrak b^*$ invariant.

Following the orbit method, in order to obtain the induced representation from the Borel subgroup $B$ to $G$, we do describe the (co)adjoint orbits of $B$ in $\mathfrak b^*$.

\begin{lem}[Orbital picture of $\mathfrak b^*$ in the nonsplit case]
The space $\mathfrak b^*$ is divided into a dijoint union of the following coadjoint $B$-orbits:
\begin{enumerate}
\item[a.] The two half-spaces $\Omega_\pm$ consisting of functionals $F = tT^* +  xX^* + yY^* +zZ^*$, 
$$\Omega_+= \{(t,x,y,z) \in \mathbb R^4 |  z >  0\}$$ 
$$\Omega_-= \{(t,x,y,z) \in \mathbb R^4 |  z <  0\}$$ 
\item[b.] A family of cylinders with hyperbolic base
$$\Omega_\alpha = \{(t,x,y,0) \in \mathbb R^3 \times \{0\} | xy = \alpha  \}, \alpha > 0$$ 
\item[c.] Four half-planes corresponding to the case $xy = 0$ but $x\ne y$
$$\Omega_{x>0} = \{t,x,0,0) \in \mathbb R^4 | x > 0\}$$
$$\Omega_{x<0} = \{t,x,0,0) \in \mathbb R^4 | x < 0\}$$
$$\Omega_{y>0} = \{t,0,y,0) \in \mathbb R^4 | y> 0\}$$
$$\Omega_{y<0} = \{t,0,y,0) \in \mathbb R^4 | y< 0\}$$
\item[d.]
The origin $$\Omega = \{(0,0,0,0) \}$$
\end{enumerate}
\end{lem}
\textsc{Proof.}
This proposition is proven by a direct computation of the (co-)adjoint action.
\hfill$\Box$

Let us now use the orbit method  \cite{dndiep1},\cite{kirillov}
Consider the linear functionals $\pm Z^*\in \Omega_\pm \subset \mathfrak g^*$ and the corresponding coadjoint orbits $\Omega_\pm = G.(\pm Z^*).$

\begin{lem}
Subalgebras $\mathfrak l = \mathbb C(X \pm iY) \oplus \mathbb CZ \subset \mathfrak u_\mathbb C$ are the positive polarizations at $\pm Z^*\in \Omega_\pm$.
\end{lem}
\textsc{Proof.}
The Lemma is proven by a direct computation the conditions from the definition of a positive polarization.
\hfill$\Box$

\subsection{Holomorphic Induction}
Following the orbit method and the holomorphic induction, we do choose the integral functionals $\lambda$, take the corresponding orbits and then choose polarization and use the holomorphic induction.

As described above, the positive root system $\Delta^+ = \Delta^+_c \cup \Delta^+_n = \{ \alpha_{kl}, 1\leq k\ne l \leq 3, \beta, 2\beta\}=\{\alpha_{12}, \alpha_{32},\alpha_{31}\}, \rho = \alpha_{32}$.
The root spaces are $\mathfrak g_\mathbb C^{\alpha_{kl}} = \mathbb C   E_{kl}$. $\mathfrak g_\beta = \mathbb R X \oplus \mathbb R Y$ and $\mathfrak g_{2\beta} = \mathbb R Z$.  
Define $$\mathfrak p_+ = \bigoplus_{\alpha \in \Delta^+_n} \mathfrak g^\alpha = \mathfrak g^{\alpha_{32}} \oplus \mathfrak g^{\alpha_{31}}=\mathbb C E_{31}\oplus \mathbb C E_{32}$$ 
and 
$$\mathfrak p_- = \bigoplus_{\alpha \in \Delta^-_n} \mathfrak g^{\alpha} = \mathfrak g^{\alpha_{23}} \oplus \mathfrak g^{\alpha_{13}}=\mathbb C E_{13}\oplus \mathbb C E_{23}$$ 

Denote $\mathcal F\subset (i\mathfrak h)^*$ the set of all linear functional $\lambda$ on $\mathfrak h_\mathbb C$ such that $(\lambda + \rho)(H_\alpha)$ is integral for any root $\alpha\in\Delta$, where $H_\alpha$ is the coroot corresponding to root $\alpha$ and $\rho$ is the half-sum of the positive roots. 
Denote also
$$\mathfrak F' = \{\lambda \in \mathfrak F | \lambda(H_\alpha) \ne 0, \forall \alpha\in \Delta \},$$
$$\mathfrak F'_0 = \{\lambda\in \mathfrak F' | \lambda(H_\alpha) > 0, \forall \alpha\in \Delta^+ _c\}$$ 
$$=\left\{\lambda\in \mathfrak F' \left| \begin{array}{l} \lambda(H_{12}) \in \mathbb N^+ \mbox{ and } \lambda(H_{31}) \in \mathbb N^+  \mbox{ (holomorphic case) or }\\
\lambda(H_{12})\in \mathbb N^+ \mbox{ and } \lambda(H_{23}) \in \mathbb N^+  \mbox{(anti-holomorphic case)  or }\\
\lambda(H_{12}) \in \mathbb N^+ \mbox{ and } \lambda(H_{13}) \in \mathbb N^+, \lambda(H_{12}) > \lambda(H_{13}) \mbox{ neither-nor case }
\end{array}
 \right .\right\}$$
 
Choose complex subalgebra $$\mathfrak e = \mathfrak p_+ \oplus \mathfrak k_\mathbb C,$$ we have
$$\mathfrak e + \overline{\mathfrak e} = \mathfrak g_\mathbb C, \mathfrak e \cap \overline{\mathfrak e} = \mathfrak k_\mathbb C $$
and therefore we have a positive polarization.

The compact root Weyl group $W_K = \langle s_{\alpha_{12}} \rangle$ is generate by a single reflection $s_{\alpha_{12}}$ then
for any $\lambda \in i\mathfrak h$, $-s_{\alpha_{12}}\lambda - \alpha_{32} = -s_{\alpha_{12}}(\lambda + \alpha_{31})$, theerefore if $V_\lambda$ is a $K$-module of lowest weight $\lambda + \rho$ then its contragradient $K$-module $V^*_\lambda$ is of heighest weight $\lambda + \alpha_{31}$.

Because $G= BK = B_1K$, the relative cohomology of  $(\mathfrak g, K)$-module with coefficients in the representation $V_\lambda$ can be reduced to the one of $B$ or $B_1=AU \subset B$ with Lie algebra $\mathfrak b_1 = \langle S= E_{13}+E_{31},X,Y,Z\rangle$. 
\begin{prop} The $(\mathfrak g, K)$-module $\pi_\lambda$ with coefficients in the representation $V_\lambda$
$$H^{q_\lambda}(G,K;\mathfrak e,V_\lambda) = H^{q_\lambda}(B, M;\mathfrak e \cap \mathfrak b,V_\lambda) =H^{q_\lambda}(B_1;\mathfrak e \cap \mathfrak b_1,V_\lambda)$$
\end{prop}

\subsection{Hochschild-Serre spectral sequence}
Remark that because in general $\mathfrak p$ is not a subalgebra, we can modify it by taking subalgebra $\mathfrak h_+ = \mathbb C(Y+iX) \oplus \mathbb C(S-iZ/2)$:
$$\mathfrak e = \mathfrak p_+ \oplus \mathfrak k_\mathbb C = \mathfrak h_+ \oplus \mathfrak k_\mathbb C, \quad\mathfrak h_+.$$ Therefore, one has
$$\mathfrak e \cap \mathfrak b_1 = \mathfrak h_+, \quad \mathfrak e \cap \mathfrak b = \mathfrak h_+ \oplus \mathfrak m_\mathbb C.$$
We may construct a Hochschild-Serre spectral sequence for this filtration.

Consider a highest weight $\lambda + \alpha_{31}$ representation  $V^*_\lambda$ of $\mathfrak k_\mathbb C$, which is trivially on $\mathfrak p_+$ extended to a representation $\xi$ of $\mathfrak e = \mathfrak p_+ \oplus \mathfrak k_\mathbb C$. The action of $\mathfrak h_+$ in $V^{\lambda + \alpha_{31}}$ is $\xi + \frac{1}{2}\tr\ad_{\mathfrak b_1}$.
Denote by $\mathcal H_\pm$ the space of representations $T_\pm$ of $B$ $\Omega_\pm$ above and by $\mathcal H_\pm^\infty$ the subspaces of smooth vectors.
Because $\dim_\mathbb C(\mathfrak p_\mathbb C) = 2$, we have $\wedge^q(\mathfrak h_+)= 0$, for all $q \geq 3$.
It is natural to define the Hochschild-Serre cobound operators
$$(\delta_\pm)_{\lambda,q} : \wedge^q(\mathfrak h_+)^* \otimes V^{\lambda + \alpha_{31}} \otimes \mathcal H_\pm^\infty \to  \wedge^{q+1}(\mathfrak h_+)^* \otimes V^{\lambda + \alpha_{31}} \otimes \mathcal H_\pm^\infty$$ and by duality their formal adjoint operators $(\delta_\pm)_{\lambda,q}^*$.
The Hochschild-Serre spectral sequence is convergent
$$\bigoplus_{r+s=q} H^r(\mathfrak e_1; H^s(M;V^{\lambda + \alpha_{31}} \otimes \mathcal H_\pm^\infty)) \Longrightarrow H^q(B;\mathfrak b_1,V_\lambda)$$

\begin{thm}
The trace of the discrete series representations in the degenerate case is a finite sum of relative traces, i. e. if $$f = \sum_{i=1}^N f_i h_i,\quad f_i \in C^\infty_c(P/U), h_i \in C^\infty_c(MA)$$ and the Hochschild-Serre spectral converges
$$\bigoplus_{p+q = n} H^1(MA; H^q(U; V))\Longrightarrow H^n(P; V)$$
$$tr \pi^\pm_n(f) = \sum_{i=1}^N \tr\sigma^\pm_k(h_i)|_{H^p(MA;\mathbb C)}  \tr\chi^\pm_k|_{ H^q(U; V))}$$
\end{thm}
\beginpf
The theorem is a consequence of the above spectral approximation.
\endpf

\section{Trace Formula}
In this section we make precise the Arthur-Selberg trace formula for $\SL(3,\bbR)$. We refer the readers to the prominent work of J. Arthur \cite{arthur}.

\subsection{Characters of unitary representations}
Let us remind that $\Gamma \subset \SL(3,\mathbb R)$ is a finitely generated  discrete subgroup with finite number of cusps,
of finite co-volume $$\vol(\Gamma\backslash \SL(3,\bbR) < \infty.$$
Let $f\in C^\infty_c(\SL(3,\bbR))$ be a smooth function of compact support. If $\varphi$ is a function from the representation space, the action of the induced representation $\inf_P^G\chi$ is the restriction of the right regular representation $R$ on the inducing space of induced representation.
$$\tr R(f)\varphi = \int_G (f(y)R(y)\varphi(x)dy) = \int_G f(y)\varphi(xy)dy$$
$$= \int_G f(x^{-1}y)\varphi(y)dy (\mbox{right invariance of Haar measure } dy) $$
$$= \int_{\Gamma\backslash G}\left(\sum_{\gamma\in \Gamma}f(x^{-1}\gamma y) 
\right)\varphi(y)dy$$

 Therefore, this action can be represented by an operator with kernel $K(x,y)$ of form
$$[R(f)\varphi](x) = \int_{\Gamma\backslash G} K_f(x,y)\varphi(y)dy,$$ where
$$K_f(x,y) = \sum_{\gamma\in \Gamma}f(x^{-1}\gamma y).$$ Because the function $f$ is of compact support, this sum convergent, and indeed is a finite sum, for any fixed $x$ and $y$ and is of class $L^2(\Gamma\backslash G \times\Gamma\backslash G )$. The operator is of trace class and it is well-known that
$$\tr R(f) = \int_{\Gamma\backslash G}K_f(x,x)dx.$$
As supposed, the discrete subgroup $\Gamma$ is finitely generated. Denote by $\{\Gamma\}$ the set of representatives of conjugacy classes. For any $\gamma\in\Gamma$ denote the centralizer of $\gamma\in \Omega \subset G$ by $\Omega_\gamma$, in particular, $G_\gamma \subset G$. Following the Fubini theorem for the doble integral, we can change the order of integration to have
$$\tr R(f) = \int_{\Gamma\backslash G} K_f(x,x)dx =\int_{\Gamma\backslash G} \sum_{\gamma\in \Gamma} f(x^{-1}\gamma x)dx$$
$$=\int_{\Gamma\backslash G} \sum_{\gamma\in \{\Gamma\}}\sum_{\delta\in \Gamma_\gamma\backslash \Gamma} f(x^{-1}\delta^{-1}\gamma\delta x)dx$$
$$=\sum_{\gamma\in \{\Gamma\}}\int_{\Gamma_\gamma\backslash G}  f(x^{-1}\gamma x)dx=\sum_{\gamma\in \{\Gamma\}}\int_{G_\gamma\backslash G}\int_{\Gamma_\gamma\backslash G_\gamma}  f(x^{-1}u^{-1}\gamma u x)dudx$$
$$= \sum_{\gamma\in \{\Gamma\}}\int_{G_\gamma\backslash G} \vol(\Gamma_\gamma \backslash G_\gamma)f(x^{-1}\gamma x)dx.$$
Therefore, in order to compute the trace formula, one needs to do:
\begin{itemize}
\item classfiy the conjugacy classes of all $\gamma$ in $\Gamma$: they are of type elliptic (different eigenvalues of the same sign),
 hyperbolic (nondegenerate, with eigenvalues of different sign),
 parabolic (denegerate)  
\item Compute the volume of form; it is the volume of the quotient of the stabilazer of the adjoint orbits. 
$\vol(\Gamma_\gamma \backslash G_\gamma)$
\item and compute the orbital integrals of form
$$\mathcal{O}(f) = \int_{G_\gamma\backslash G} f(x^{-1}\gamma x)d\dot x$$ 
The idea is to reduce these integrals to smaller endoscopic subgroups in order to the correponding integrals are ordinary or almost  ordinary. 
\end{itemize}

\subsection{Stable trace formula}
The main difficult in the previous section is that the sum of traces which are unstable under the action of the Galois group.
We refer the readers to \cite{arthur} for more details.

\section{Endoscopy}
In this section, the method that was used by J.-P. Labesse for $\SL(2,\mathbb R)$ is applied to analyze our case of $\SL(3,\mathbb R)$. The result is in the same way obtained, cf. \cite{labesse}.
\subsection{Orbital integrals}
\textsl{The simplest case} is the elliptic case when $\gamma = \diag(a_1,a_2, a_3), \quad a_1a_2a_3 = 1$ and they are pairwise different. In this case, because of Iwasawa decomposition $x=mauk$, and the $K$-bivariance,  the orbital integral is
$$\mathcal O_\gamma(f) = \int_{G_\gamma\backslash G} f(x^{-1}\gamma x)dx = \int_U f(u^{-1}\gamma u)du=$$ $$\int_{\mathbb R^3} f(\begin{pmatrix} 1 & x & z\\ 0 & 1& y \\ 0 & 0 &  1\end{pmatrix}^{-1}\begin{pmatrix} a_1 & 0 & 0\\ 0 & a_2 & 0 \\ 0 & 0 & a_3 \end{pmatrix} \begin{pmatrix} 1 & x & z\\ 0 & 1& y \\ 0 & 0 &  1\end{pmatrix})dxdydz=$$ 
$$\int_{\mathbb R^3} f(\begin{pmatrix} 1 & -x & yx - z\\ 0 & 1 & -y\\ 0 & 0 & 1\end{pmatrix}\begin{pmatrix} a_1 & 0 & 0\\ 0 & a_2 & 0\\ 0 & 0 & a_3 \end{pmatrix} \begin{pmatrix} 1 & x & z\\ 0 & 1& y\\ 0 & 0 & 1\end{pmatrix})dxdydz$$ $$\int_{\mathbb R^3} f(\begin{pmatrix} 1 & -x & yx - z\\ 0 & 1 & -y\\ 0 & 0 & 1\end{pmatrix}\begin{pmatrix} a_1  &a_1x & a_1z\\ 0 &a_2 & a_2y\\ 0 & 0 & a_3\end{pmatrix})dxdydz = $$ 
$$\int_{\mathbb R^3} f(\begin{pmatrix} a_1  &(a_1-a_2)x &(\dots )x + (\dots) + (a_1-a_3)z y \\ 0 &a_2 & (a_2-a_3)y\\ 0 & 0 & a_3\end{pmatrix})dxdydz = $$ 
$$|a_1-a_2|^{-1}|a_2 -a_3|^{-1}|a_1 - a_3|^{-1}  f^H.$$
The integral is absolutely and uniformly convergent and therefore is smooth function of $a\in (\bbR^*_+)^2$. Therefore the function $$f^H(\gamma) = \Delta(\gamma)^{-1}\mathcal O_\gamma(f), \quad \Delta(\gamma) =\prod_{1\leq i<j\leq 3} |a_i- a_j|$$ is a smooth function on the endoscopic group $H= (\bbR^*)^2$. 

\textsl{The second case} is the case where $\gamma = \begin{pmatrix} k_\theta & 0\\ 0 & 1\end{pmatrix} = \begin{pmatrix} \cos\theta & \sin\theta & 0\\ -\sin\theta & \cos\theta & 0\\ 0 & 0 &1 \end{pmatrix}$. We have again, $x= mauk$. $a =\diag(a_1,a_2,a_3), a_1a_2a_3 = 1$ and
$$\mathcal O_{k(\theta)}(f) = \int_{G_{k(\theta)}\backslash G} f(k^{-1}u^{-1}a^{-1}m^{-1}k(\theta)mauk)dmdudadk=$$
$$ \int_{G_{k(\theta)}\backslash G} f(u^{-1}a^{-1}m^{-1}k(\theta)mau)dmduda=$$
$$  \int_{G_{k(\theta)}\backslash G} f(\begin{pmatrix} 1 & -x &yx-z\\ 0 & 1 & -y\\ 0 & 0 & 1 \end{pmatrix}\begin{pmatrix} a_1^{-1} & 0 & 0\\ 0 & a_2^{-1} & 0\\ 0 & 0 & a_3^{-1} \end{pmatrix} \begin{pmatrix} \cos\theta & \sin\theta & 0\\ -\sin\theta & \cos\theta & 0\\ 0 & 0 & 1 \end{pmatrix}$$ $$\times \begin{pmatrix} a_1 &0 & 0\\ 0 & a_2 & 0 \\ 0 & 0 & a_3 \end{pmatrix}\begin{pmatrix} 1 & x & z \\ 0 & 1 & y \\ 0 & 0 & 1  \end{pmatrix})duda dk(\theta)=$$
$$  \int_{G_{k(\theta)}\backslash G} f(\begin{pmatrix} 1 & -x &yx-z\\ 0 & 1 & -y\\ 0 & 0 & 1 \end{pmatrix} \begin{pmatrix} \cos\theta & a_2a_1^{-1}\sin\theta & 0\\ -a_2^{-1}a_1\sin\theta & \cos\theta & 0\\ 0 & 0 & 1 \end{pmatrix} \begin{pmatrix} 1 & x & z \\ 0 & 1 & y \\ 0 & 0 & 1  \end{pmatrix})duda dk(\theta)=$$
$$ c\int_1^\infty f(\begin{pmatrix}\cos\theta & t_1\sin\theta & 0\\ -t_1^{-1}\sin\theta & \cos\theta & 0\\ 0 & 0 & 1\end{pmatrix})\prod_{i=1}^2|t_i-t_i^{-1}|\frac{dt_i}{t_i}=$$ 
$$c \int_0^{+\infty}\sgn(t-1) \tilde{f}(\begin{pmatrix}\cos\theta & t\sin\theta & 0\\ -t^{-1}\sin\theta & \cos\theta & 0\\ 0 & 0 & 1 \end{pmatrix})dt,$$ 
where $c$ is some constant and $\tilde{f}$ some good function
When $f$ is an element of the Hecke algebra, i.e. $f$ is of class $C^\infty_0(G)$ and is $K$-bivariant, the integral is converging absolutely and uniformly. Therefore the result is a function $F(\sin\theta)$. 
The function $f$ has compact support, then the integral is well convergent at $+\infty$. At the another point $0$, we develope the function $F$ into the Tayor-Lagrange of the first order with respect to $\lambda = \sin\theta \to 0$
$$F(\lambda) = A(\lambda) + \lambda B(\lambda),$$ where $A(\lambda) = F(0)$ and $B(\lambda)$ is the error-correction term $F'(\tau)$ at some intermediate value $\tau, 0 \leq \tau \leq t$.
Remark that 
$$ \begin{pmatrix} \sqrt{1-\lambda^2} & t\lambda & 0\\ -t^{-1}\lambda & \sqrt{1-\lambda^2}& 0 \\ 0 & 0 & 1 \end{pmatrix}$$ 
$$=\begin{pmatrix} t^{1/2}&0& 0\\ 0 & t^{-1/2} & 0\\
 0 & 0 & 1 \end{pmatrix} \begin{pmatrix} \sqrt{1-\lambda^2} & \lambda & 0\\ \lambda & \sqrt{1-\lambda^2} & 0\\ 0 & 0 & 1 \end{pmatrix}\begin{pmatrix} t^{-1/2} & 0 & 0\\ 0 &t^{1/2} & 0\\
0 & 0 & 1 \end{pmatrix} $$ we have
$$B= \frac{dF(\tau)}{d\lambda} =\frac{d}{d\lambda} \left. \int_0^{+\infty} \sgn(t-1)f(\begin{pmatrix} \sqrt{1-\lambda^2} & t\lambda & 0\\ -t^{-1}\lambda & \sqrt{1-\lambda^2} & 0\\  0 & 0 & 1 \end{pmatrix})dt \right|_{t=\tau}$$
$$=  \int_0^{+\infty} \sgn(t-1)g(\begin{pmatrix} \sqrt{1-\lambda^2} & t\lambda & 0\\ -t^{-1}\lambda & \sqrt{1-\lambda^2}& 0\\ 0 & 0 & 1 \end{pmatrix})\frac{dt}{t}, $$
where $g\in C^\infty_c(N)$
and $g(\lambda) \cong O(-t^{-1}\lambda)^{-1}.$ $B$ is of logarithmic growth and $$B(\lambda) \cong \ln(|\lambda|^{-1})g(1)$$ up to constant term,  and therefore is contimuous. 
$$A = F(0) = |\lambda|^{-1} \int_0^\infty f(\begin{pmatrix}1 & \sgn(\lambda)u & 0\\ 0 & 1 & 0\\ 0 & 0 & 1\end{pmatrix}) du - 2f(I_3) + o(\lambda)$$
Hence the functions 
$$G(\lambda) = |\lambda|(F(\lambda) + F(\lambda)),$$
$$H(\lambda) = \lambda(F(\lambda) - F(-\lambda))$$ have the Fourier decomposition
$$G(\lambda) = \sum_{n=0}^N (a_n|\lambda|^{-1} + b_n)\lambda^{2n} + o(\lambda^{2N})$$
$$H(\lambda) = \sum_{n=0}^N h_n\lambda^{2n} + o(\lambda^{2N})$$
Summarizing the discussion, we have that in the case of $\gamma = k(\theta)$, there exists also a continuous function $f^H$ such that $$f^H(\gamma) = \Delta(\gamma) (\mathcal O_\gamma(f) - \mathcal O_{w\gamma}(f))=\Delta(k(\theta))\mathcal {SO}_\gamma(f), $$ where $\Delta(k(\theta)) = -2i\sin\theta$.

\begin{thm}\cite{labesse}
There is a natural function $\varepsilon : \Pi \to \pm 1$ such that in the Grothendieck group of discrete series representation ring, $$\sigma_G = \sum_{\pi\in \Pi} \varepsilon(\pi)\pi,$$ the map $\sigma \mapsto \sigma_G$ is dual to the map of geometric transfer, that for any $f$  on $G$, there is a unique $f^H$ on $H$
$$\tr \sigma_G(f) = \tr\sigma(f^H).$$
\end{thm}
\textsc{Proof.}
There is a natural bijection $\Pi_\mu \cong \mathfrak D(\mathbb R,H,G)$, we get a pairing
$$\langle .,.\rangle : \Pi_\mu \times \mathfrak k(\mathbb R,H,G)\to \mathbb C.$$ Therefore we have
$$\tr \Sigma_\nu(f^H) = \sum_{\pi\in \Pi_\Sigma} \langle s,\pi\rangle \tr \pi(f).$$
\hfill$\Box$

Suppose given a complete set of endoscopic groups $H = \mathbb S^1 \times \mathbb S^1 \times \{\pm 1\}$ or $\SL(2,\mathbb R) \times \{\pm 1\}$.  For each group, there is a natural inclusion
$$\eta: {}^LH \hookrightarrow {}^LG$$

Let $\varphi: DW_\bbR  \to {}^LG$ be the Langlands parameter, i.e. a homomorphism from the Weil-Deligne group $DW_\bbR = W_\bbR\ltimes \bbR^*_+$ the Langlands dual group,
$\bbS_\varphi$ be  the set of conjugacy classes of Langlands parameters modulo the connected component of identity map. For any $s\in 
\bbS_\varphi$,
$\check{H}_s = \Cent(s,\check G)^\circ$ the connected component of the centralizer of $s\in\bbS_\varphi$ we have $\check{H}_s$ is conjugate with $H$. 
Following the D. Shelstad pairing $$\langle s, \pi\rangle : \bbS_\varphi \times \Pi(\varphi) \to \bbC$$
$$\varepsilon(\pi) = c(s)\langle s,\pi\rangle.$$
Therefore, the relation
$$\sum_{\sigma\in \Sigma_s} \tr \sigma(f^H) = \sum_{\pi\in\Pi} \varepsilon(\pi) \tr \pi(f)$$
can be rewritten as
$$\widetilde{\Sigma}_s(f^H) = \sum_{s\in\Pi} \langle s,\pi\rangle \tr\pi(f)$$
and
$$\widetilde{\Sigma}_s(f^H) = c(s)^{-1}\sum_{\sigma\in\widetilde{\Sigma}_s} \tr \sigma(f^H) .$$
We arrive, finally to the result
\begin{thm}\cite{labesse}
$$\tr \pi(f) = \frac{1}{\#\bbS_\varphi}\sum_{s\in\bbS_\varphi} \langle s,\pi\rangle\widetilde{\Sigma}_s(\check{f}^H).$$
\end{thm}

In the Langlands picture of the trace formula, the trace of the restriction of the regular representation on the cuspidal parabolic part is the coincidence of of the spectral side and the geometric side.
\begin{equation}
\sum_{\pi} m(\pi) \hat{f}(\pi) = \sum_{\gamma\in \Gamma\cap H} a^G_\gamma \hat{f}(\gamma) \end{equation}
Let us do this in  more details.

\subsection{Stable orbital integrals}

Let us remind that the \textit{orbital integral} is defined as
$$\mathcal{O}(f) = \int_{G_\gamma\backslash G} f(x^{-1}\gamma x)d\dot x$$

The complex Weyl group is isomorphic to $\mathfrak S_3$ while the real Weyl group
is isomorphic to $\mathfrak S_2$ . The set of conjugacy classes inside a strongly regular
stable elliptic conjugacy class is in bijection with the pointed set
$\mathfrak S_3 /\mathfrak S_2$
 that can be viewed as a sub-pointed-set of the group
$\mathfrak E(\mathbb R, T, G) = (Z_2 )^2$
We shall denote by $\mathfrak K(\mathbb R, T, G)$ its Pontryagin dual.

Consider $\kappa \ne 1$ in $\mathfrak K(\mathbb R, T, G)$ such that $\kappa(H_{13}) = -1$. 
 Such a $\kappa$ is unique:
in fact one has necessarily
$\kappa(H_{12}) = \kappa(H_{13}) = -1$.

The endoscopic group $H$ one associates to $\kappa$ is isomorphic to $\SL(2,\bbR)$ and

can be embedded in $G$ as 
$$\begin{pmatrix}
ua & iub & 0\\  -iuc & ud & 0\\  0 & 0 & 1 \end{pmatrix},$$ $$ w = \begin{pmatrix} a & b\\ c & d \end{pmatrix}, ad-bc = 1 \mbox{ and } u = \pm 1.$$
%It will be useful to consider also the twofold cover
%$H_1 = S(U (1) \times U (1)) \times \SL(2)$.

Let $f_\mu$ be a pseudo-coefficient for the discrete series representation $\pi_\mu$ then the \textit{$\kappa$-orbital integral} of a  regular element $\gamma$
in $T (R)$ is given by
$$\mathcal O^\kappa_\gamma(f_\mu) =\int_{G_\gamma\backslash G} \kappa(x) f_\mu(x^{-1}\gamma x) d\dot x $$ 
$$= \sum_{\sign(w) =1} \kappa(w) \Theta^G_\mu(\gamma^{-1}_w)  = \sum_{\sign(w) =1} \kappa(w)\Theta_{w\mu}(\gamma^{-1}),$$ because there is a natural bijection between the left coset classes and the right coset classes.

\subsection{Endoscopic transfer}
The transfer factor $\Delta(\gamma, \gamma_H )$ is given by
$$\Delta(\gamma, \gamma_H )=(−1)^{q(G)+q(H)} \chi_{G,H} (\gamma)\Delta_B (\gamma^{−1} ) . \Delta_{B_H}(\gamma_H^{-1})^{-1}$$
for some character $\chi_{G,H}$ defined as follows. Let $\xi$ be a character of the twofold covering $\mathfrak h_1$ of $\mathfrak h$, then
$\chi_{G,H} (\gamma^{−1} ) = e^{\gamma^{\rho - \rho_H +\xi}}$  defines a character of $H$, corresponding to $\mathfrak h$, because it is trivial on any fiber of the cover.

With such a choice we get when $\sign (w) = 1$
and $w\ne 1$, we have $\kappa(w) = -1$ and 
$$\Delta(\gamma^{-1},\gamma_H^{-1})\Theta_{w\mu}^G(\gamma) = -\frac{\gamma_H^{w\mu + \xi} - \gamma_H^{w_0w\mu + \xi}}{\gamma^{\rho_H}\Delta_{B_H}(\gamma_H}$$ therefore
$$\Delta(\gamma,\gamma_H)\Theta^G_{w\mu} (\gamma^{-1}) = \kappa(w)^{-1}\mathcal{SO}_\nu^H(\gamma_H^{-1}),$$ where $\nu = w\mu +\xi$ is running over the corresponding $L$-package of discrete series representations for the endoscopic group $H$.
Therefore we have the following formula
$$\Delta(\gamma,\gamma_H){\mathcal O}_\gamma^\kappa(f_\mu) = \sum_{\nu=w\mu + \xi\atop \sign(w) = 1} {\mathcal {SO}}_\nu^H(\gamma_H^{-1})$$ or
$$\Delta(\gamma,\gamma_H){\mathcal O}_\gamma^\kappa(f_\mu) = \sum_{\nu=w\mu + \xi\atop \sign(w) = 1} {\mathcal {SO}}_{\gamma_H}(g_\nu), $$ where $g_\nu$ is pseudo-coefficient for any one of the discrete series representation of the endoscopic subgroup $H$ in the $L$-package of $\mu$.

For any $$f^H = \sum_{\nu=w\mu+\rho\atop \sign(w) = 1} a(w,\nu)g_\nu,\quad a(w_1,w_2\mu = \kappa(w_2) \kappa(w_2w_1)^{-1}$$ we have the formula
$$\tr \Sigma_\nu (f^H) = \sum_{w} a(w,\nu) \tr \pi_{w\mu}(f).$$

\section{Poisson Summation Formula}

\subsection{Endoscopic orbital integrals}

\begin{thm}\cite{labesse}
There is a function $\varepsilon: \Pi \to \pm 1$
such that, if we consider $σ_G$ in the Grothendieck group defined by
$$\sigma_G =\sum_{\pi \in\Pi}\varepsilon(\pi)\pi, $$
then $\sigma \mapsto \sigma_G$ is the dual of the geometric transfer:
$$\tr \sigma_G (f ) = \tr \sigma(f^H )$$
\end{thm}
\textsc{Proof.}
There is a natural bijection $\Pi_\mu \cong \mathfrak D(\mathbb R,H,G)$, we get a pairing
$$\langle .,.\rangle : \Pi_\mu \times \mathfrak k(\mathbb R,H,G)\to \mathbb C$$. Therefore we have
$$\tr \Sigma_\nu(f^H) = \sum_{\pi\in \Pi_\Sigma} \langle s,\pi\rangle \tr \pi(f).$$
\hfill$\Box$

\subsection{Endoscopic Trace Formula}
\begin{thm}
$$\tr R(f)_{L^2_{cusp}(\Gamma \backslash \SL(3,\bbR)} = \sum_{\Pi_\mu}\sum_{\pi\in\Pi_\mu}m(\pi)\mathcal S\Theta_\pi(f) = \sum_{ \Pi_\mu} \Delta(\gamma,\gamma_H)\mathcal {SO}(f_\mu),$$ where
$$\mathcal S\Theta_\pi(f) = \sum_{\pi\in \Pi}\kappa(\pi) \Theta_\pi(f)$$ is the sum of Harish-Chandra characters of the discrete series running over the stable conjugacy classes of $\pi$ 
and $$\mathcal {SO}(f_\mu) = \sum_{\lambda\in \Pi_\mu}\kappa(\pi_\lambda) \mathcal O(f_\lambda)$$ is the sum of orbital integrals weighted by a character $\kappa : \Pi_\mu \to \{\pm1\}$. 
\end{thm}
\beginpf\;
The proof just is  a combination of the previous theorems.
\endpf

%${}^1$ \textsc{Institute of Mathematics, Vietnam Academy of Sciences and Technology, 18 Hoang Quoc Viet Road, Cau Giay District, 10307 Hanoi, Vietnam}\\
%\texttt{Email: dndiep@math.ac.vn}\\
%\\
%${}^2$ \textsc{Medicine College, Thai Nguyen Universtity, Thai Nguyen City, Vietnam}\\
%\texttt{Email: phuongquynhtn@gmail.com}

\begin{thebibliography}{AAAA}
\bibitem[A]{arthur}{\sc J. Arthur}, {\it An Introduction to the Trace Formula}, Clay Mathematics Proceedings "Harmonic Analysis,
THe Trace Formula, and Shimura Varieties", Volume 4, 2005, pp. 01-259.

\bibitem[D1]{dndiep1}{\sc Do Ngoc Diep}, {\it Methods of Noncommutative Geometry for Group C*-Algebras }, Chapman \& Hall/ Research Notes in Mathematics Series, Vol. 416, 346 pp., Boca Raton, Florida, New York, London, Berlin, 1999.

\bibitem[F]{fomin}\textsc{A. I. Fomin}, \textit{Semisimple irreducible representations of $\SL(3,R)$}, Funkt. Anal. ego  Prilozh., 9(1975), No 3, 67-74.

\bibitem[K]{kirillov}{\sc A.A. Kirillov}, {\it Elements of the Theory of Representations}, Springer Verlag, Berlin - Heidelberg - New York, 1975.

\bibitem[La]{labesse} \textsc{J.-P. Labesse}, \textit{Itroduction to Endoscopy},  Snowbird Lectures, June 2006
www.institut.math.jussieu.fr/projets/fa/bpFiles/Labesse.pdf‎

\bibitem[L1]{langlands1}{\sc R. Langlands}, {\it Orbital Integrals on Forms of SL(3), I}, sunsite.ubc.ca/DigitalMathArchive/Langlands/orbintI/orb-ps.ps;  Amer. Jour. Math. 105 (1983), 465-506.

\bibitem[L2]{langlands2}{\sc R.P. Langlands, D. Shelstad}, {\it Orbital Integrals on Forms of SL(3) , II }, sunsite.ubc.ca/DigitalMathArchive/Langlands/pdf/orb2-ps.pdf


\end{thebibliography}
\end{document}